\documentclass[a4paper]{article}

\usepackage[utf8]{inputenc}
\usepackage{lmodern}
\usepackage{amssymb,amsfonts,amsthm,amsmath,bbm,mathrsfs,aliascnt,mathtools}
\usepackage[english]{babel}
\usepackage[colorlinks]{hyperref}
\usepackage{tikz}
\usetikzlibrary{decorations}
\usetikzlibrary{decorations.pathreplacing}
\tikzset{individu/.style={draw,thick}}

\theoremstyle{plain}
\newtheorem{theorem}{Theorem}[section]
\newtheorem{corollary}[theorem]{Corollary}
\newtheorem{lemma}[theorem]{Lemma}
\newtheorem{proposition}[theorem]{Proposition}

\theoremstyle{definition}

\theoremstyle{remark}
\newtheorem{remark}[theorem]{Remark}

\numberwithin{equation}{section}

\usepackage{tikz}

\makeatletter \newcommand \listoftodos{\section*{Todo list} \@starttoc{tdo}}
\newcommand\l@todo[2]
{\par\noindent \textit{#2}, \parbox{10cm}{#1}\par} \makeatother

\newcommand{\R}{\mathbb{R}}

\newcommand{\ind}[1]{\mathbbm{1}_{\left\{#1\right\}}}

\renewcommand{\bar}[1]{\overline{#1}}
\renewcommand{\tilde}[1]{\widetilde{#1}}
\renewcommand{\hat}[1]{\widehat{#1}}

\newcommand{\e}{\mathrm{e}}
\newcommand{\dd}{\mathrm{d}}

\DeclareMathOperator{\E}{\mathbb{E}}
\newcommand{\Q}{\mathbb{Q}}
\renewcommand{\P}{\mathbb{P}}

\renewcommand{\epsilon}{\varepsilon}

\title{On the local times of \\ noise reinforced Bessel processes}
\author{Jean Bertoin\thanks{Institute of Mathematics, University of Zurich, Switzerland.} }
\date{
\small{Dedicated to the memory of Marc Yor}
}

\begin{document}

\maketitle

\begin{abstract} We investigate the effects of noise reinforcement on a Bessel process of dimension $d\in(0,2)$, and more specifically on the asymptotic behavior of its additive functionals. This leads us to introduce a local time process and its  inverse. We identify the latter as an increasing self-similar (time-homogeneous) Markov process, and from this, several explicit results can be deduced.
\end{abstract}


\medskip


\section{Introduction}
\label{sec:introduction}
Loosely speaking, the notion of stochastic reinforcement refers to certain modifications of the dynamics of a given random process, such that, depending on whether the reinforcement is positive or negative, some of its previous steps are either more likely or less likely to be repeated in the future. See the survey \cite{Pem} and references therein for background. In the  continuous time setting, the name noise reinforcement is meant to suggest that reinforcement acts on the increments of the process over infinitesimal time intervals. 

A noise reinforced Brownian motion is a process   $\hat B=(\hat B_t)_{t\geq 0}$ which can be defined as 
$$\hat B_t = \frac{t^{p}}{\sqrt{1-2p}}{B_{t^{1-2p}}},$$
where $B=(B_t)_{t\geq 0}$ is  a standard Brownian motion  and $p\in(-\infty,1/2)$ the reinforcement parameter. It has  notably appeared in scaling limit theorems for two-color urn models \cite{Gouet, BHZ}, for certain simple  random walks with memory called elephant random walks \cite{BaurBer, ColGavSch2}, and more recently, as the universal weak limit of step-reinforced random walks (see \cite{BerVla} for $p\in(0,1/2)$ and \cite{BerRos} for $p\in(-1,1/2)$). We also mention that  $\hat B$  is a building block of so-called  noise-reinforced L\'evy processes \cite{NRLP}. The restriction for the range of $p$ can be explained informally by the fact that  for $p\geq 1/2$, the reinforcement would be too strong and induce an instantaneous  explosion.

We define analogously a  \emph{noise reinforced Bessel process} (of dimension $d>0$ and with reinforcement parameter $p\in(-\infty,1/2)$),   $\hat R=(\hat R_t)_{t\geq 0}$,  by
\begin{equation} \label{E:NRBES}
\hat R_t = \frac{t^{p}}{\sqrt{1-2p}}{R_{t^{1-2p}}},
\end{equation}
where $R=(R_t)_{t\geq 0}$ denotes a Bessel process of dimension $d$ started from $R_0=0$; we refer to  \cite[Section XI.1]{RY} and \cite{Lawler} for background. 
We agree implicitly in \eqref{E:NRBES} that $\hat R_0=0$ even when $p<0$, as this should be plain from the H\"older property of the sample paths of Bessel processes. 
 When $d$ is an integer, $\hat R$ can thus be viewed as the Euclidean norm of a $d$-dimensional noise reinforced Brownian motion. Observe that for any $d>0$, $\hat R$ is a  self-similar  time-inhomogeneous Markov process. Although this will not be needed in the sequel,  the reader may find usefull to recall that when $d>1$, the Bessel process solves the stochastic differential equation
$$\dd R_t= \dd B_t + \frac{d-1}{2R_t}  \dd t,$$
 where $B$ denotes  some Brownian motion. It follows readily from \eqref{E:NRBES} that, in turn, the noise reinforced Bessel process is a solution to
 $$\dd \hat R_t=\dd B'_t +  \left( \frac{d-1}{2 \hat R_t} + p \frac{ \hat R_t}{t} \right) \dd t,$$
 where $B'$ is another Brownian motion. Hence, noise reinforcement amounts to adding  a drift term
 given by the product of the time-average velocity $\hat R_t/t$ and the reinforcement parameter $p$ to the dynamics of the Bessel process.
 When $0< d < 1$, the Bessel process is no longer a semi-martingale, and the stochastic differential equations above would need to be interpreted properly; see \cite{BerBes} and \cite[Chapter 10]{MansuyYor}.

Loosely speaking, an important issue in the literature on reinforced stochastic processes is to understand of how the reinforcement impacts ergodic properties. 
The analysis of the asymptotic behavior of additive functionals  is a classical theme for null
recurrent time-homogeneous Markov processes, see notably \cite{KR,DK,PSV,KK} and \cite[Section XIII.2]{RY} 
(see also  \cite{Kono} and \cite{HNX} for fractional Brownian motions), and the purpose of the present work is to study this question for noise reinforced Bessel processes. We focus on dimensions $d\in(0,2)$, as Bessel processes are transient for dimensions $d>2$, and, by the construction \eqref{E:NRBES}, transience is easily transferred to their noise reinforced versions.

We assume from now on that $0<d<2$ and will use 
$$\alpha= 1-d/2 \in(0,1)$$
as a more convenient parameter.

\begin{theorem}\label{T1} There exists  a  continuous non-decreasing process $\hat{L} = (\hat{L}_t)_{t\geq 0}$,  called the local time at level $0$ of $\hat R$, which  satisfies 
$$\mathrm{Supp}\, \dd \hat{L}_t =\{t\geq0: \hat R_t=0\}, \qquad \text{almost surely,}$$
and such that the following assertions hold for every function  $f: \R_+\to \R$  in $ \mathrm{L}^1(x^{1-2\alpha} \dd x)$.

\begin{enumerate}
\item[(i)]
There is the convergence in distribution on $\mathcal{C}(\R_+,\R)$ 
as $n\to \infty$:
$$\left(n^{-\alpha} \int_0^{nt} f(\hat R_s) \dd s\right)_{t\geq 0} \Longrightarrow  c_1(f)\hat{L},$$
where
$$c_1(f)= \alpha^{-1} \int_0^{\infty} f(x) x^{1-2\alpha} \dd y.$$

\item[(ii)]
Furthermore, if $c_1(f)=0$ and $f$ has compact support, then
$$\left(n^{-\alpha/2} \int_0^{nt} f(\hat R_s) \dd s\right)_{t\geq 0} \Longrightarrow  \left(\beta_{c_2(f)\hat{L}(t)}\right)_{t\geq 0},$$
where
$$c_2(f)= \frac{4}{\alpha} \int_0^{\infty} \left(\int_{0} ^{x} f(t) t^{1-2\alpha} \dd t \right )^2 x^{2\alpha-1}\dd x.$$
and $\beta = (\beta_t)_{t\geq 0}$ is a Brownian motion independent of $\hat{L}$.
\end{enumerate}
\end{theorem}
Theorem \ref{T1} may be a bit disappointing as the reinforcement parameter $p$ does not seem to play much of a role. In particular, noise reinforcement does not affect the growth exponent $\alpha$ of additive functionals, which might nonetheless come as a surprise for the following reason.
For time-homogeneous point-recurrent Markov processes, the growth exponent of additive functionals
is related to the decay of the tail distribution of return times, and one can readily check from \eqref{E:NRBES} that  the tail distribution of the duration of the excursion of $\hat R$, say straddling time $1$, decays with exponent $-\alpha(1-2p)$, a quantity which does depend on the reinforcement parameter $p$.

Actually, the impact of noise reinforcement is perceived at the level of  the distribution of the local time $\hat{L}$. 
\begin{theorem}\label{T2} 
The one-dimensional distributions  of $\hat{L}$ are determined by their entire moments which are given by 
$$\E(\hat{L}_t)=\frac{(2t)^{\alpha } }{\Gamma(1-\alpha)  (1-2p)^{\alpha -1}}$$
and,  for every  $n\geq 2$,
$$\E(\hat{L}_t^n)= (2t)^{\alpha n} (1-2p)^{1-\alpha n}
 \frac{\Gamma(1+\alpha)^{n-1}}{ \Gamma(1-\alpha)^n} \Gamma(n)  
 \prod_{k=1}^{n-1}  \frac{\Gamma(\alpha k/(1-2p))}{\Gamma( \alpha(1+k/(1-2p)))}.
 $$
 \end{theorem}
The formula of Theorem \ref{T2} becomes much simpler  when there is no reinforcement,  that is for $p=0$. Indeed it then reduces to
$$\left( \frac{2^{\alpha}\Gamma(1+\alpha)}{\Gamma(1-\alpha)}\right)^n \frac{n!}{\Gamma(1+\alpha n)},$$
and we recognize the $n$-th moment of the local time (at level $0$) of the Bessel process $R$ as we should expect; see \cite[Proposition 3.4]{DRVY}. For $p\neq 0$,  $\hat{L}_1$ is not proportional to a Mittag-Leffler variable and  Theorems \ref{T1} and \ref{T2} hence contrast with the Darling-Kac theorem \cite{DK} for time-homogeneous Markov processes. 

How to properly  define the local time $\hat{L}$ for the noise reinforced processes clearly lies the heart of this work.
This classical question is well understood for time-homogeneous Markov processes:  the local time  is a continuous additive functional of the Markov process that increases only when the latter visits a given state (which must be regular for itself), and this characterizes the local time at this state up to a constant factor; see  \cite[Definition V.3.12 and Theorem VI.3.13]{BlGe}. The question may however be ambiguous for time-inhomogeneous Markov processes.
In our setting, a hasty guess would be to argue that, as by \eqref{E:NRBES} the zero set of $\hat R$ 
is the image of the zero set of $R$ by the map $t\mapsto t^{1-2p}$, one might simply define
the local time $\hat{L}_t$ of $\hat R$ as the time change ${L}_{t^{1-2p}}$ of the local time ${L}$ of $R$. 
However this construction overlooks the distortion of space induced by the dilation with factor $t^p/\sqrt{1-2p}$ in \eqref{E:NRBES}, and does not yield the correct notion, at least for the problems that we are interested in here. 

The local time $\hat{L}$ of the reinforced Bessel process will be constructed in Section \ref{sec:RLT} via a noise reinforced version of the celebrated Tanaka formula. This will enable us
to establish Theorem \ref{T1} in Section \ref{sec:ProofThm1} by classic stochastic calculations. In Section \ref{sec:inverseRLT}, we will study the inverse local time process $\hat \lambda$; the key observation is that $\hat \lambda$ is an increasing self-similar time-homogeneous Markov process (recall that the inverse of the local time of a time-homogeneous Markov process is always a subordinator). The study of self-similar time-homogeneous Markov processes on $\R_+$ has been initiated by Lamperti \cite{Lamperti} and this topic has been very active during the last three decades or so. The proof of Theorem \ref{T2} will be given in Section \ref{sec:inverseRLT}, using some nowadays well known results in this field. 
Finally,  in Section \ref{sec:occ},  we will construct  the two-parameter family $(\hat{L}^x_t)$ of local times of $\hat R$ (that is at time $t$ and level $x$) as occupation densities, and study their regularity.

Last but not least, most of the ideas in this text stem from works of Marc Yor and his co-authors, as the reader will quickly realize by browsing through the rest of this paper or having a glance at the references. Let me seize this occasion to acknowledge once more the great influence that Marc has had on my own research.

\section{Local time after noise reinforcement}
\label{sec:RLT}
We start by observing  from the construction \eqref{E:NRBES} and the scaling property of Bessel processes that for every $s> 0$, $\hat R_s$ has the same law as $R_{s/(1-2p)}$. The latter is absolutely continuous with density   
\begin{equation} \label{E:Besdens}
\frac{2^{\alpha}}{ \Gamma(1-\alpha)}\left( \frac{s}{1-2p}\right)^{\alpha-1}
x^{1-2\alpha} \exp\left (-\frac{(1-2p)x^2}{2s}\right ), \qquad x>0;
\end{equation}
 see \cite[page 446]{RY}. In particular, $\hat R_1^q\in \mathrm{L}^1(\P)$  for all $q>2\alpha-2$, and in that case,
 \begin{equation}\label{E:momBes}
 \E(\hat R_s^q)=s^{q/2} \E(\hat R_1^q)<\infty \qquad \text{for all }s> 0.
\end{equation}

We next lift 
 from \cite{DRVY} some features of ${L}=({L}_t)_{t\geq 0}$,  the local time  at level $0$ of the Bessel process. We stress that the latter is taken in the sense of Markov processes (when $d>1$, even though $R$ is a semimartingale, its  local time in the sense of semimartingales  \cite[Chapter VI]{RY} is identically zero). Plainly, the zero-set of $R$ has zero Lebesgue a.s. and the Stieltjes measure $\dd {L}_t$ is almost surely singular with respect to the Lebesgue measure on $\R_+$. 
According to \cite[Theorem 2.1]{DRVY}, ${L}$  is conveniently characterized by an analog of Tanaka's formula. Specifically, the process $R^{2 \alpha}$ is a nonnegative submartingale with Doob-Meyer decomposition 
\begin{equation}\label{E:Tanaka}
R^{2 \alpha}_t =  2\alpha \int_0^t R^{2\alpha-1}_s \dd B_s+ {L}_t,
\end{equation}
where, as usual, $B$ denotes a standard Brownian motion. We also recall from \cite[Proposition 3.1]{DRVY} that
\begin{equation} \label{E:TLO}
{L}_t=\lim_{\varepsilon \to 0+} \frac{2\alpha(1-\alpha)}{\varepsilon^{2-2\alpha} } \int_0^t \ind{R_s\leq \varepsilon} \dd s.
\end{equation}

The noise reinforced Bessel process $\hat R$ is only a time-inhomogeneous Markov process and  thus  there is \textit{a priori} no `canonical' way of defining its local time at the level $0$. Nonetheless, if we think of $\hat R$ as a `perturbation' of the Bessel process $R$, \eqref{E:Tanaka} points at the following.

\begin{lemma}[Noise reinforced Tanaka's formula] \label{L1} The process $\hat R^{2\alpha}$ is a continuous semimartingale with canonical decomposition 
$$\hat R_t^{2\alpha}= 2\alpha \int_0^{t} \hat R^{2\alpha-1}_s \dd B'_s + \hat V_t,$$
where $B'$ is a Brownian motion  and $\hat V$ a process with bounded variation. More precisely,
 the canonical decomposition of $\hat V$ as the sum of its absolutely continuous and its singular components (with respect to the Lebesgue measure on $\R_+$) is given by
$$\hat V_t =  2\alpha p \int_0^{t} \frac{\hat R^{2\alpha}_s}{s} \dd s + \hat{L}_t,$$
where 
\begin{equation}\label{E:RLT}
\hat{L}_t= (1-2p)^{-\alpha} \int_0^{t^{1-2p}} s^{2\alpha p/(1-2p)}  \dd {L}_s. 
\end{equation}
We henceforth refer to $\hat{L}=(\hat{L}_t)_{t\geq 0}$ as the local time (at level $0$) of $\hat R$. 
\end{lemma}
\begin{remark} The expression \eqref{E:RLT} has appeared recently in the case $\alpha = 1/2$ (i.e. $d=1$) as scaling limit for the number of zeros of an elephant random walk  in diffusive regimes, see \cite[Theorem 3.1]{Zero}. Informally, the scaling limit of an elephant random walk is a noise reinforced Brownian motion $\hat B$, and therefore the scaling limit of its number of zeros should be the local time of $\hat B$ at $0$. 
\end{remark}

\begin{proof}
The calculations below are related to the proof of  \cite[Proposition 2.1]{HY}, see also \cite[Exercise 6.11]{CY}. 
We first write
$$(1-2p)^{\alpha} \hat R^{2\alpha}_t =   \left( t^{1-2p}\right )^{2\alpha p/(1-2p)} R^{2\alpha}_{t^{1-2p}}$$
and apply It\^{o}'s formula on the right-hand side (use \eqref{E:Tanaka} and beware of the singularity at the initial time when $p<0$)  to get for every $0<\varepsilon<t$:
\begin{align*}
(1-2p)^{\alpha} \left( \hat R^{2\alpha}_t -  \hat R^{2\alpha}_{\varepsilon} \right)
= \, &   2\alpha \int_{\varepsilon^{1-2p}}^{t^{1-2p}} s^{2\alpha p/(1-2p)} R^{2\alpha-1}_s \dd B_s
 \\ +&  \int_{\varepsilon^{1-2p}}^{t^{1-2p}} s^{2\alpha p/(1-2p)}  \dd {L}_s \\
+& \frac {2\alpha p}{1-2p}  \int_{\varepsilon^{1-2p}}^{t^{1-2p}} s^{2\alpha p/(1-2p)} R^{2\alpha}_s \frac{\dd s}{s}. 
\end{align*}

We next  observe that the change of variables $s=r^{1-2p}$ yields 
$$\frac {2\alpha p}{(1-2p)^{\alpha+1}}  \int_{\varepsilon^{1-2p}}^{t^{1-2p}} s^{2\alpha p/(1-2p)} R^{2\alpha}_s \frac{\dd s}{s}= 2\alpha p \int_{\varepsilon}^{t} \frac{\hat R^{2\alpha}_r}{r} \dd r.$$
Recall that the sample paths of the Bessel process $R$ are a.s. H\"{o}lder-continuous with exponent $1/2-\eta$ for any $\eta>0$, so by \eqref{E:NRBES} the same holds for $\hat R$. This enables us to take the limit as $\epsilon\to 0+$ in the preceding integral, and we get the continuous process
$$2\alpha p \int_{0}^{t} \frac{\hat R^{2\alpha}_r}{r} \dd r, \qquad t\geq 0.$$

Similarly, we have
$$\frac {2\alpha }{(1-2p)^{\alpha}} \int_{\varepsilon^{1-2p}}^{t^{1-2p}} s^{2\alpha p/(1-2p)} R^{2\alpha-1}_s \dd B_s= 2\alpha \int_{\varepsilon}^{t} \hat R^{2\alpha-1}_r \dd B'_r,$$
where
\begin{equation}\label{E:B'}
B'_t= (1-2p)^{-1/2} \int_0^{t^{1-2p}} s^{p/(1-2p)} \dd B_{s}
\end{equation}
is another Brownian motion. Using \eqref{E:momBes}, we can take the limit as $\varepsilon\to 0+$, which yields the continuous martingale
$$2\alpha \int_{0}^{t} \hat R^{2\alpha-1}_r \dd B'_r,\qquad t\geq 0.$$

Since obviously $ \hat R^{2\alpha}_{\varepsilon}$ converges to $0$ a.s. as $\varepsilon \to 0+$, we now see that \eqref{E:RLT} defines a continuous non-decreasing process,
and its Stieltjes measure $\dd \hat{L}_t$ is singular with respect to the Lebesgue measure since the same holds for $\dd {L}_t$. Furthermore, it should be plain from \eqref{E:NRBES} that the support of $\dd \hat{L}_t$ coincides  with the zero set of $\hat R$ a.s.
\end{proof} 

It is seen from the reinforced Tanaka formula that $\hat{L}$ is self-similar with exponent $\alpha$. Recall further that the Bessel local time ${L}$ has finite moments of any order which can be computed explicitly, see \cite[Proposition 3.4]{DRVY}. An integration by parts  in   \eqref{E:RLT} gives
$$(1-2p)^{\alpha} \hat{L}_t= t^{2\alpha p} {L}_{t^{1-2p}} + 2\alpha p \int _0^{t} {L}_{s^{1-2p}} s^{2\alpha p-1} \dd s,$$
and it follows that $\hat{L}_1\in {\mathrm L}^q(\P)$ for any $q>0$ as well.
For future use, we record that, by self-similarity, 
\begin{equation} \label{E:momL}
\E(\hat{L}_t^q) =  t^{q \alpha} \E(\hat{L}_1^q)< \infty\qquad \text{for all }t\geq 0.
\end{equation}
At this point the determination of $\E(\hat{L}_1^q)$, even for integer values of $q$, does not seem easy.

It follows also readily from Lemma \ref{L1}, \eqref{E:momBes}, \eqref{E:momL} and the Burkholder-Davis-Gundy  inequalities, that if we set
$\hat R^*_t=\sup_{0\leq s \leq t} \hat R_s$, then $\hat R^*_t\in {\mathrm L}^q(\P)$ for all $q\geq 0$. By self-similarity, we thus have
\begin{equation} \label{E:momR}
\E((\hat R^*_t)^q) =  t^{q/2} \E((\hat R^*_1)^q) < \infty  \qquad \text{for all }t\geq 0.\end{equation}
We shall make use of these observations in the next section. 

\section{Proof of Theorem \ref{T1}}
\label{sec:ProofThm1}
In this section, we will derive Theorem \ref{T1} from the noise reinforced Tanaka formula and self-similarity. The argument is standard in stochastic calculus and some easy details of the calculation will be left to the reader.

\begin{lemma} \label{L2} Let $g: \R_+\to \R$ in $\mathrm{L}^1(\dd x)$, write $\bar g(0)=\int_0^{\infty}g(x)\dd x$. The following convergence holds in $\mathrm{L}^1(\P)$.
$$\lim_{T\to \infty} T^{-\alpha} \sup_{0\leq t \leq T} \left| 2\alpha^2 \int_0^{t} g(\hat R^{2\alpha}_s)\hat R_s^{4\alpha-2} \dd s - \bar g(0) \hat{L}_t  \right| = 0 . $$
\end{lemma}
\begin{proof} We can  assume that $g\geq 0$ without loss of generality. Define for $x\geq 0$
$$\bar g(x)=\int_x^{\infty} g(y)\dd y\quad \text{and} \quad G(x) = \int_0^x \bar g(y) \dd y.$$
Thus $G$ is a concave non-decreasing function with first derivative $G'=\bar g$ and second derivative  in the $\mathrm{L}^1(\dd x)$-sense $G'' =-g$. An application of It\^{o}'s formula to  
$G(\hat R_t^{2\alpha})$ using the semimartingale decomposition of $\hat R^{2\alpha}$
in Lemma \ref{L1} yields
\begin{align} \label{E:decomp1}& 2\alpha^2 \int_0^t g(\hat R^{2\alpha}_s)\hat R_s^{4\alpha-2} \dd s  
- \bar g(0) \hat{L}_t  \nonumber \\ 
= &\, - G(\hat R_t^{2\alpha}) + 2\alpha \int_0^t \bar g(\hat R^{2\alpha}_s) \hat R_s^{2\alpha-1} \dd B'_s 
+ 2\alpha \int_0^t \bar g(\hat R^{2\alpha}_s)  \frac{\hat R^{2\alpha}_s}{s} \dd s.
\end{align}

First,  the identity $\sup_{t\leq T} G(\hat R_t^{2\alpha})= G((\hat R^*_t)^{2\alpha})$,  where $ \hat R^*$ denotes the running supremum process of $\hat R$, \eqref{E:momR}  and the fact that $G(x)=o(x)$ as $x\to \infty$ 
yield
$$\lim_{T\to \infty} \E\left( T^{-\alpha} \sup_{t\leq T} G(\hat R_t^{2\alpha})\right)=0.$$

Next, fix $\varepsilon>0$ and choose $a>0$ sufficiently large so that $\bar g(a)\leq \varepsilon$. 
Denote the stochastic integral  in the right-hand side of \eqref{E:decomp1} by $M_t$.
Its quadratic variation can be bounded from above by 
$$\langle M \rangle_t\leq  \bar g(0)^2\int_0^t  \mathbbm{1}_{\hat R^{2\alpha}_s\leq a} \hat R_s^{4\alpha-2} \dd s 
+ \varepsilon^2   \int_0^t  \hat R_s^{4\alpha-2} \dd s.$$ 
 On the one hand, we readily deduce from the explicit form of the density \eqref{E:Besdens} that 
 $$ \E\left(  \mathbbm{1}_{\hat R^{2\alpha}_s\leq  a} \hat R_s^{4\alpha-2}\right)= O(s^{\alpha-1}) \qquad \text{as }s\to \infty.$$
On the other hand, \eqref{E:momBes} entails that the expectation of the second integral above satisfies 
$$ \E\left( \int_0^t  \hat R_s^{4\alpha-2} \dd s\right) = O(t^{2\alpha})\qquad \text{as }t\to \infty.$$
Since $\varepsilon$ can be arbitrarily small, 
this shows that
$$\lim_{T\to \infty} \E( T^{-2\alpha} \langle M \rangle_T)=0.$$
We conclude from the Burkholder-Davis-Gundy inequality that
 $$\lim_{T\to \infty} \E\left( T^{-\alpha} \sup_{t\leq T} |M_t|\right)=0.$$

For the third term in the right-hand side of  \eqref{E:decomp1}, let $\varepsilon$ and $a$ be as above. 
We then split the integral to get the bound
$$\int_0^T \bar g(\hat R^{2\alpha}_s) \hat R^{2\alpha}_s \frac{\dd s}{s}
\leq \bar g(0) \int_0^1  \hat R^{2\alpha}_s \frac{\dd s}{s}
+\bar g(0) a \int_1^T  \frac{\dd s}{s} + \varepsilon  \int_1^T  \hat R^{2\alpha}_s \frac{\dd s}{s},
$$
from which we immediately infer as above  that 
$$\lim_{T\to \infty}  \E\left(T^{-\alpha} \int_0^T \bar g(\hat R^{2\alpha}_s) \hat R^{2\alpha}_s \frac{\dd s}{s}\right)=0.$$
This completes the proof. 
\end{proof}

The first  claim of Theorem \ref{T1} can now be seen from the self-similarity of $\hat{L}$ by an application of Lemma \ref{L2} to the function $g(y)=f(y^{1/2\alpha}) y^{-2+1/\alpha}$, so that
$f(r)= g(r^{2\alpha}) r^{4\alpha-2}$ and 
$$\bar g(0)=\int_0^{\infty} g(y) \dd y = 2\alpha \int_0^{\infty} f(x) x^{1-2\alpha} \dd x = 2\alpha^2 c_1(f).$$
 In order to establish the second part of the statement, we need the following refinement of Lemma \ref{L2} when $\bar g(0)=0$. 

\begin{lemma} \label{L3} Let $g: \R_+\to \R$ in $\mathrm{L}^1(\dd x)$ with compact support, write $\bar g(x)=\int_x^{\infty}g(y)\dd y$. If $\bar g(0)=0$, then the following convergence holds in $\mathrm{L}^1(\P)$.
$$\lim_{T\to \infty} T^{-\alpha/2} \sup_{0\leq t \leq T} \left| \alpha \int_0^{t} g(\hat R^{2\alpha}_s)\hat R_s^{4\alpha-2} \dd s - \int_0^t \bar g(\hat R^{2\alpha}_s) \hat R_s^{2\alpha-1} \dd B'_s \right| = 0 . $$
\end{lemma}
\begin{proof} We use the same notation as in the proof of Lemma \ref{L2}. As $\bar g(0)=0$, 
\eqref{E:decomp1} now reads
\begin{align*} & \int_0^t \bar g(\hat R^{2\alpha}_s) \hat R_s^{2\alpha-1} \dd B'_s-\alpha \int_0^{t} g(\hat R^{2\alpha}_s)\hat R_s^{4\alpha-2} \dd s \nonumber \\
= &\,  \frac{1}{2\alpha} G(\hat R_t^{2\alpha}) 
-  \int_0^t \bar g(\hat R^{2\alpha}_s) \hat R^{2\alpha}_s \frac{\dd s}{s}.
\end{align*}
We have obviously $$\lim_{T\to \infty} \E\left( T^{-\alpha/2} \sup_{t\leq T} G(\hat R_t^{2\alpha})\right)=0,$$
as $G$ is now bounded.
Moreover the function $x \mapsto x^{2\alpha}\bar g(x^{2\alpha})$ is also bounded, and it follows that
$$\int_0^T \E(|\bar g(\hat R^{2\alpha}_s)| \hat R^{2\alpha}_s)\frac{\dd s}{s} = O(\log T) \qquad \text{as }T\to \infty.$$
This yields the claim. 
\end{proof}
We also need the following weak limit theorem for the stochastic integral  in Lemma \ref{L3}.

\begin{lemma} \label{L4} With the same assumptions and notation as in Lemma \ref{L3},  there
is the weak convergence in $\mathcal{C}(\R_+,\R)$ as $n\to \infty$
$$\left(n^{-\alpha/2} \int_0^{nt}  \bar g(\hat R^{2\alpha}_s) \hat R_s^{2\alpha-1} \dd B'_s \right)_{t\geq 0} \Longrightarrow  \left(\beta_{c\hat{L}(t)}\right)_{t\geq 0},$$
where
$$c= \frac{1}{2 \alpha^2}\int_0^{\infty} \bar g(x)^2 \dd x$$
and $\beta = (\beta_t)_{t\geq 0}$ is a Brownian motion independent of $\hat{L}$.
\end{lemma}
\begin{proof} The argument is again standard; see \cite[Section XIII.2]{RY}. To start with, we observe that
the reinforced Bessel process $\hat R$, and hence also its local time $\hat{L}$, are measurable with respect to the Brownian motion $B'$. Indeed, recall the Tanaka formula \eqref{E:Tanaka} for the Bessel process $R$. Writing $R^2=(R^{2\alpha})^{1/\alpha}$ and applying It\^{o}'s formula, we see that $R$ is measurable with respect to the Brownian motion $B$; see \cite[page 439]{RY}. Our assertion is now plain from \eqref{E:NRBES} and \eqref{E:B'}. 

Next, for every $n\geq 1$, we write 
$\beta^n$ for the Dambis-Dubins-Schwarz  Brownian motion associated to the stochastic integral of the statement. Calculations similar to those already performed in this section for quadratic variations enable us to apply \cite[Theorem XIII.2.3]{RY} and  show that the sequence $ \left ( n^{-1/2} B'_{nt}, \beta^n_t\right)_{t\geq 0}$ converges in distribution to a pair of independent Brownian motions, say $(\beta',\beta)$. We can then conclude as in the proof of \cite[Theorem XIII.2.6]{RY}, using the fact that, thanks to Lemma \ref{L2}, the sequence of rescaled Brownian motions and quadratic variation processes 
$$\left( n^{-1/2} B'_{nt}, n^{-\alpha} \int_0^{nt}  \bar g(\hat R^{2\alpha}_s)^2 \hat R_s^{4\alpha-2} \dd s\right)_{t\geq 0} $$
converges in law as $n\to \infty$ towards $(B',c \hat{L})$. 
\end{proof}
The second claim of Theorem \ref{T1} is now plain from  an application of Lemmas \ref{L3} and  \ref{L4} to the function $g(y)=f(y^{1/2\alpha}) y^{-2+1/\alpha}$, so that
$$\bar g(x)=\int_x^{\infty} g(y) \dd y =\int_0^x g(y) \dd y = 2\alpha \int_0^{x^{1/2\alpha}} f(t) t^{1-2\alpha} \dd t,$$
and then 
$$\frac{1}{2 \alpha^4}\int_0^{\infty} \bar g(x)^2 \dd x=4\alpha^{-1} \int_0^{\infty} \left(\int_{0} ^{x} f(t) t^{1-2\alpha} \dd t \right )^2 x^{2\alpha-1}\dd x
= c_2(f).$$

\section{Inverse of noise reinforced local time as a self-similar Markov process}
\label{sec:inverseRLT}
The main purpose of this section is to establish Theorem \ref{T2}. By self-similarity, it suffices to verify the formula for the entire moments there  for $t=1$. Our approach will rely on properties of increasing self-similar  Markov processes. This  looks rather indirect, and one might expect that Kac's moment formula (see \cite{FP}) should provide a simpler way. Unfortunately, I have not been able to make direct calculations  explicit, likely by clumsiness. On the other hand, the present approach sheds light on the fine structure  of $\hat{L}$ and may be interesting on its own right. 

From place to place, it will be convenient to write $X(t)$ instead of $X_t$, where $X$ stands for some stochastic process depending on the time parameter $t$.
We introduce the right-continuous inverse  process of the noise reinforced local time
$$\hat \lambda_t=\inf\{s\geq 0: \hat{L}_s>s\}, \qquad t\geq 0.$$
Similarly, we write 
$$ \lambda_t=\inf\{s\geq 0:  {L}_s>s\}, \qquad t\geq 0,$$
for the inverse of the Bessel local time 
and recall from \cite[Proposition 3.2]{DRVY} that the latter is an $\alpha$-stable subordinator. 
 The construction \eqref{E:RLT} of $\hat{L}$ in terms of ${L}$ readily yields a simple expression for $\hat \lambda$ in terms of $\lambda$. 
In this direction, recall that an  $\alpha$-stable subordinator grows roughly like a power function of time with exponent $1/\alpha$, 
and therefore, informally speaking,  
$$\lambda_s^{2\alpha p/(1-2p)}\approx s^{2p/(1-2p)}.$$ Since $2p/(1-2p)>-1$,  the map
$$ t\mapsto (1-2p)^{-\alpha}  \int_0^{t}  \lambda_s^{2\alpha p/(1-2p)} \dd s, \qquad t\geq 0,$$
is bijective on $\R_+$ a.s. Its inverse process $\tau$ is defined  implicitly by
\begin{equation}\label{E:rho}
 \int_0^{\tau(t)}  \lambda_s^{2\alpha p/(1-2p)}  \dd s = (1-2p)^{\alpha}  t ,  \qquad t\geq 0.
 \end{equation}

\begin{lemma} \label{L5} With probability one, there is the identity
$$\hat \lambda_t= \lambda_{\tau(t)}^{1/(1-2p)}\qquad \text{for all } t\geq 0.$$

\end{lemma}
\begin{proof} 
Since $L\circ \lambda = \mathrm{Id}$,   \eqref{E:RLT} can be rewritten as
$$
(1-2p)^{\alpha} \hat{L}_t=   \int_0^{L(t^{1-2p})}  \lambda_s^{2\alpha p/(1-2p)} \dd s, 
$$
from which we deduce
$$L\left(\hat \lambda_t^{1-2p}\right)=\tau(t).$$
Hence 
$$ \lambda_{\tau(t)-}^{1/(1-2p)}\leq \hat \lambda_t \leq  \lambda_{\tau(t)}^{1/(1-2p)}, \qquad \text{for all } t\geq 0,$$
and since $\hat \lambda$ is right-continuous, we obtain the formula of the statement. 
 \end{proof}
 
 Lemma \ref{L5} enables us to identify the distribution of the inverse local time $\hat \lambda$ as a self-similar Markov process;
 in turn this will be our main tool for establishing Theorem \ref{T2}. 
 
 \begin{corollary}\label{C1} 
 \begin{enumerate}
 \item[(i)]
 The inverse   $\hat \lambda$ of the noise reinforced local time is an increasing self-similar Feller process on $[0,\infty)$
 that starts from the boundary point $0$ and has scaling exponent $\alpha$. 
 
 \item[(ii)] The infinitesimal generator of $\hat \lambda$ is given on $(0,\infty)$ by
 $$ \hat{\mathcal G}f(x)  = (1-2p)^{\alpha+1} \frac{2^{-\alpha}}{\Gamma(\alpha)}  x^{-\alpha}  \int_1^{\infty} \left( f(x v)-f(x)\right)  v^{-2p} (v^{1-2p}-1)^{-\alpha - 1} \dd v,
$$
with $f:[0,\infty)\to \R$ a generic function in $\mathcal{C}^1_0$ (i.e. $f$ is  continuously differentiable on $[0,\infty)$, and both $f$ and $f'$ vanish at infinity). 
  \end{enumerate}
 \end{corollary}
 \begin{proof} (i)
The stable subordinator $\lambda$ is a Feller process on $\R_+$, and the time-change $\tau$ has been defined as the inverse of a (perfect continuous homogeneous) additive functional. According to  e.g. \cite[Section III.21]{WillRo},  $\lambda\circ \tau$ is strongly Markovian, and the same holds for
$\hat \lambda=(\lambda\circ \tau)^{1/(1-2p)}$ since the map $x\mapsto x^{1/(1-2p)}$ is bijective on $\R_+$. 

On the other hand, recall  that $\hat{L}$ inherits self-similarity from ${L}$. This entails that the right-continuous inverse $\hat \lambda$ is in turn self-similar with scaling exponent $\alpha$. 
The Feller property of $\hat \lambda$ follows (see \cite[Theorem 2.1]{Lamperti} on $(0,\infty)$ and \cite[Theorem 1]{BY} including the boundary point $0$).

(ii)  The infinitesimal generator $\hat{\mathcal G}$ of $\hat \lambda$ is computed accordingly.
 First, recall from \cite[Proposition 3.2]{DRVY} that the stable subordinator $\lambda$ has no drift and L\'evy measure 
 $2^{-\alpha} \Gamma(\alpha)^{-1} t^{-\alpha - 1} \dd t $ on $(0,\infty)$. Hence that its infinitesimal generator ${\mathcal G}$ is characterized  on $(0,\infty)$ by
 \begin{align*}
 {\mathcal G}f(x) &= \frac{2^{-\alpha}}{\Gamma(\alpha)}\int_0^{\infty} \left( f(x+t)-f(x)\right)  t^{-\alpha - 1} \dd t \\
 &= \frac{(2x)^{-\alpha}}{\Gamma(\alpha)}\int_1^{\infty} \left( f(ux)-f(x)\right)  (u-1)^{-\alpha - 1} \dd u \\
\end{align*}
for any $f\in \mathcal{C}^1_0$; see \cite[Theorem 31.5 on page 208]{Sato}. 
 According to Volkonskii's formula \cite[III.21.4 on page 277]{WillRo} and \eqref{E:rho}, the infinitesimal generator $\tilde {\mathcal G}$ of the time-changed process $\tilde \lambda=\lambda \circ \tau$
 is
 $$
 \tilde {\mathcal G}f(x) = (1-2p)^{\alpha} \frac{2^{-\alpha}}{\Gamma(\alpha)}  x^{-\alpha /(1-2p)}  \int_1^{\infty} \left( f(ux)-f(x)\right)  (u-1)^{-\alpha - 1} \dd u.
$$
Since, according to Lemma \ref{L3}, $\hat \lambda =\tilde \lambda^{1/(1-2p)}$, we readily conclude that
$$
 \hat{\mathcal G}f(x) = (1-2p)^{\alpha} \frac{2^{-\alpha}}{\Gamma(\alpha)}  x^{-\alpha}  \int_1^{\infty} \left( f(x u^{1/(1-2p)})-f(x)\right)  (u-1)^{-\alpha - 1} \dd u.
 $$
 The change of variables $u=v^{1-2p}$ completes the proof of the statement.
    \end{proof}

    As a classical consequence of  self-similarity, 
    \begin{equation}\label{E:samelaw} 
    \hat{L}_1 \text{ and }\hat \lambda_1^{-\alpha} \text{ have the same distribution,}
    \end{equation}
     and in particular, we have for every positive integer $n$:
  $$\E( \hat{L}_1^n)= \E( \hat \lambda_1^{-\alpha n}).$$
 As  we know from Corollary \ref{C1}(i) that  $\hat \lambda$ is an increasing self-similar Markov process started from $0$,  the right-hand side  can be computed from \cite[Theorem 1]{BeCa}. 
 Specifically, one finds 
 $$\E( \hat \lambda_1^{-\alpha n}) = \left\{
 \begin{matrix}
 1/(\alpha m) & \text{ for $n=1$}\\
  \Gamma(n)/(\alpha m \hat \Phi(\alpha) \cdots \hat \Phi(\alpha(n-1))) & \text{for }n\geq 2,
 \end{matrix} 
 \right.$$
  where $\hat \Phi$ stands for the Laplace exponent of the subordinator $\hat \xi$ associate to the increasing self-similar Markov process $\hat \lambda$ by the Lamperti transformation (see \cite[Chapter 5]{KP} for background), and
  $m=\hat \Phi'(0)$. The proof of Theorem \ref{T2} thus readily reduces  to verifying the identity
  \begin{equation} \label{E:hatphi}
  \hat \Phi(r) =   2^{-\alpha} (1-2p)^{\alpha}  \frac{\Gamma(1-\alpha)}{ \alpha \mathrm B(\alpha,r/(1-2p))}  ,  
  \end{equation}
 where $\mathrm{B}(\cdot, \cdot)$ denotes the beta function, since then 
  $$ \alpha m= \alpha \hat \Phi'(0) =  2^{-\alpha} (1-2p)^{\alpha}   \frac{ \Gamma(1-\alpha)}{ 1-2p}.$$

 \begin{proof}[Proof of \eqref{E:hatphi}] 
 According to \cite[Proposition 3.2]{DRVY}, the process $(\lambda_{at})_{t\geq 0}$,
 with 
 $$a=2^{\alpha}\Gamma(\alpha+1)/\Gamma(1-\alpha),$$
 is a standard $\alpha$-stable subordinator, that is with Laplace exponent $\kappa(r)=r^{\alpha}$.
We know from \cite[Theorem 5.11(i)]{KP} that its Lamperti subordinator is a so-called $\beta$-subordinator with Laplace exponent
$\kappa^*(r)= \Gamma(r+\alpha)/\Gamma(r)$. It follows that  the Lamperti subordinator $\xi$
of the Bessel inverse local time $\lambda$ has the Laplace exponent 
$$\Phi(r)= a^{-1} \kappa^*(r)=  2^{-\alpha}  \frac{\Gamma(1-\alpha)}{ \Gamma(1+\alpha)} \, \frac{\Gamma(r+\alpha)}{\Gamma(r)}.$$

Raising to the power $(1-2p)^{-1}$ turns $\lambda$ into $\lambda^{1/(1-2p)}$, which is again a self-similar Markov process, now with scaling exponent $\alpha(1-2p)$ and Lamperti subordinator $(1-2p)^{-1}\xi$. 
Finally, we perform the time-change given by \eqref{E:rho}. It is readily checked that
$\hat \lambda = \lambda^{1/(1-2p)}\circ \tau$ is still  a self-similar Markov process with scaling exponent $\alpha(1-2p)+ 2\alpha p = \alpha$ and Lamperti subordinator $\hat \xi$ given by
$$\hat \xi_t = (1-2p)^{-1}\xi_{(1-2p)^{\alpha}t}, \qquad t\geq 0.$$
The latter has  Laplace exponent 
$$  \hat \Phi(r)= (1-2p)^{\alpha} 
 \Phi(r/(1-2p)) =   2^{-\alpha} (1-2p)^{\alpha}  \frac{\Gamma(1-\alpha)}{ \Gamma(1+\alpha)} \, \frac{\Gamma(\alpha+r/(1-2p))}{\Gamma(r/(1-2p))} .$$ 
We obtain \eqref{E:hatphi} after simplifying the expression above using the beta function. Alternatively, one can also check the formula from Corollary \ref{C1}(ii) for the infinitesimal generator of $\hat \lambda$; see \cite{Lamperti} and \cite{CaCh}. 
   \end{proof}

We now conclude this section by presenting some further applications of self-similar Markov processes to the law of $\hat{L}_1$. Following \cite{CPY}, we introduce the so-called the exponential functional
\begin{equation} \label{E:Ihat}
\hat I= \int_0^{\infty} \exp(-\alpha \hat \xi_t) \dd t,
\end{equation}
where $\hat \xi$ denotes the subordinator with Laplace exponent 
$\hat \Phi$ given by \eqref{E:hatphi}.
\begin{corollary}\label{C2}
The law of $\hat I$ coincides with the size-biased distribution of $\hat{L}_1$, that is the identity
$$\E(f(\hat I)) = (1/2-p)^{-\alpha}   \frac{ 1-2p}{ \Gamma(1-\alpha)} \E(\hat{L}_1 f(\hat{L}_1))$$
holds for all measurable functions $f: \R_+\to \R_+$.
\end{corollary}
\begin{proof} This follows immediately from \eqref{E:samelaw} and  \cite[Theorem 1]{BY}. 
\end{proof}
Corollary \ref{C2} enables one to obtain a number of properties of the distribution of $\hat{L}_1$ from the literature on exponential functionals of L\'evy processes. In particular, it is know (see \cite{CPY, PRVS}) that the law of $I$ is absolutely continuous, and that its density 
$$k(x)=\P(\hat I\in \dd x)/\dd x, \qquad x>0,$$solves a certain equation given in terms of the tail of the L\'evy measure of $\hat \xi$.  As it is plain from \eqref{E:RLT} that the law of $\hat{L}_1$ has no atom at $0$, Corollary \ref{C2} entails that the latter is absolutely continuous with density proportional to $x^{-1}k(x)$. 
In this direction, recall from  \cite[Theorem 2.4]{PaSa} that $k$ is infinitely differentiable on $(0,\infty)$.
The study of the asymptotic behavior of $k(x)$ as $x\to \infty$  culminates with the recent manuscripts \cite{Haas} and \cite{MinSa}. In particular, one readily deduces from Corollary \ref{C1}(ii)  that the subordinator $\hat \xi$ has no drift and L\'evy measure 
$$\hat \pi(\dd x) = (1-2p)^{\alpha+1} \frac{2^{-\alpha}}{\Gamma(\alpha)}  \left(1-\e^{-(1-2p)x}\right)^{-\alpha-1} \e^{-\alpha(1-2p)x} \dd x, \qquad x>0.$$
The latter is referred to as an $(a,b,c)$-L\'evy measure in \cite{Haas}, and Lemma 21 there provides a sharp estimate of $k(x)$ in this setting.

Many more properties of $k$ are discussed in these two papers \cite{Haas,MinSa}; and the interested the reader will find further important references in their bibliographies.

\section{Occupation densities }\label{sec:occ}
The purpose of this final section is to point at  the existence of jointly continuous local times for noise reinforced Bessel processes, at least provided that one focusses on strictly positive levels.

\begin{proposition} \label{P1} There exists a two-parameter process $\left(\hat{L}^x_t\right)_{t\geq 0, x>0}$ which is jointly continuous a.s. and such that for every $t\geq0$ and  every Borel function $h: \R_+\to \R_+$, we have the occupation density formula
$$\int_0^t h(\hat R_s) \dd s = \alpha^{-1} \int_0^{\infty} h(x) \hat{L}^x_t x^{1-2\alpha} \dd x.$$
\end{proposition}
\begin{proof}
It is easy to see that for any $t>0$,  the law of the Bessel process $R$ and that of its noise reinforced version $\hat R$, both viewed as continuous processes on the time-interval $[0,t]$, are mutually singular. 
We start by observing that nonetheless,  this singularity disappears as soon as a small neighborhood of the initial time is discarded.

More precisely, fix some $\eta>0$. We readily deduce from the noise reinforced Tanaka formula (recall also that $\alpha=1-d/2$) that
$$\hat R_{t+\eta}^2= \hat R^2_{\eta} + 2\int_{\eta}^{t+\eta} \hat R_s \dd B'_s + 2p \int_{\eta}^{t+\eta} \frac{\hat R^2_s}{s} \dd s + dt. $$
Then introduce 
$$B''_s= B'_{s+\eta} - B'_{\eta} + p\int_{\eta}^{s+\eta} \ind{r>\eta} \frac{\hat R_r}{r} \dd r,\qquad s\geq 0,$$
so that 
\begin{equation} \label{E:EDS} 
\hat R_{t+\eta}^2= \hat R^2_{\eta} + 2\int_{0}^t \hat R_{s+\eta} \dd B''_s + d t. 
\end{equation}
We can construct, by a Girsanov transformation,  a law $\Q_{\eta}$ which is equivalent to the original probability law $\P$ on the sigma field $\sigma(\hat R_s, s\leq t+\eta)$, and such that $B''$ is a Brownian motion under $\Q_{\eta}$. From  \eqref{E:EDS} and the uniqueness of the solution of the stochastic differential equation for square Bessel processes (see \cite[Section XI.1]{RY}), we see that under $\Q_{\eta}$, the process 
$(\hat R_{\eta+s}^2)_{s\geq 0}$ is a square Bessel process of dimension $d$, started from 
$\hat R_{\eta}^2$. 

We then deduce from   \cite[Proposition 3.1]{DRVY} that for every $\eta>0$, there exists a two-parameter process $\left(\hat{L}^x_{\eta,t}\right)_{t\geq \eta,x\geq 0}$, which is jointly continuous  a.s. (with respect to $\Q_{\eta}$, and hence also with respect to $\P$), and such that for every Borel function $h: \R_+\to \R_+$ and $t>\eta$, we have the occupation density formula
\begin{equation} \label{E:occdens}\int_{\eta}^t h(\hat R_s) \dd s = \alpha^{-1} \int_0^{\infty} h(x) \hat{L}^x_{\eta,t} x^{1-2\alpha} \dd x.
\end{equation} 

It should be plain that the map $\eta\mapsto \hat{L}^x_{\eta,t}$ increases as $\eta$ decreases.
This enables us to define $ \hat{L}^x_{t} = \lim_{\eta \to 0+}  \hat{L}^x_{\eta,t}$ for every $x\geq 0$ and $t\geq 0$. Our claim now follows from \eqref{E:occdens}  by monotone convergence
and the easy fact that $\hat{L}^x_{\eta,t}=\hat{L}^x_{t}$ provided that
$\hat R^*_{\eta}=\sup_{0\leq s \leq \eta} \hat R_s <x$.
\end{proof}

We expect that $ \hat{L}^x_t $ converges almost surely to $ \hat{L}_t$ as $x\to 0+$, but we have only  been able to prove that the convergence holds in $\mathrm{L}^1(\P)$.
\begin{proposition}\label{P2} We have for every $t\geq 0$ that
$$
\lim_{x\to 0+} \hat{L}^x_t= \hat{L}_t\qquad \text{in }\mathrm{L}^1(\P).
$$
\end{proposition} 
\begin{proof} We use the notation in the proof of Proposition \ref{P1}, and claim first that  for every $\eta>0$, there is the identity
\begin{equation} \label{E:approxtl}
\hat{L}^0_{\eta, t}= \hat{L}_t-\hat{L}_{\eta}.
\end{equation}
To see this, consider a  function $f:\R_+\to \R_+$ in $\mathrm{L}^1(x^{1-2\alpha}\dd x)$  such that $c_1(f)=1$ (in the notation of Theorem \ref{T1}(i)). For every $n\geq 1$,
$n^{2-2\alpha} f(nx) x^{1-2\alpha} \dd x$ thus defines a probability measure on $\R_+$, and  there is the weak convergence 
 \begin{equation} \label{E:weakconv}
 \lim_{n\to \infty} n^{2-2\alpha} f(nx) x^{1-2\alpha} \dd x = \delta_0(\dd x),
 \end{equation}
 where $\delta_0$ denotes the Dirac mass at $0$. 
From \eqref{E:occdens} and the  continuity of $x\mapsto \hat L^x_{\eta,t}$ at $x=0$, we get
$$\hat{L}^0_{\eta, t}= \lim_{n\to \infty}n^{2-2\alpha}   \int_{\eta}^t f(n\hat R_s) \dd s. $$
Next, using \eqref{E:NRBES}, we have
$$\int_{\eta}^t f(n\hat R_s) \dd s 
= \frac{1}{1-2p} \int_{{\eta}^{1-2p}}^{t^{1-2p}} r^{2p/(1-2p)} f\left( \frac{n}{\sqrt{1-2p}} R_r\right) \dd r.
$$
The Bessel process $R$ possesses jointly continuous local times $(L^x_t)_{x,t\geq 0}$ in the sense of \cite[Proposition 3.1(ii)]{DRVY}. By considering distribution functions, we deduce
as above from \eqref{E:weakconv},  that with probability one, 
$$\lim_{n\to \infty}  n^{2-2\alpha} f\left( \frac{n}{\sqrt{1-2p}} R_r\right) \dd r = (1-2p)^{1-\alpha} \dd L^0_r,$$
in the sense of the vague convergence of Radon measures on $[0,\infty)$.
Since $L^0=L$  according to \cite[Proposition 3.1(i)]{DRVY}, we conclude that
$$\hat{L}^0_{\eta, t}= (1-2p)^{-\alpha} \int_{{\eta}^{1-2p}}^{t^{1-2p}} r^{2p/(1-2p)} \dd L_r = \hat{L}_t-\hat{L}_{\eta},$$
where the second equality is \eqref{E:RLT}. 

We next deduce from \eqref{E:approxtl},  the continuity of $x\mapsto \hat{L}^{x}_{\eta, t}$ at $x=0$, the inequality $\hat{L}_t^x \geq \hat{L}_{\eta,t}^{x}$, and the continuity of $\eta \mapsto \hat L_{\eta}$  at $\eta=0$,  that 
\begin{equation} \label{E:liminf} 
\liminf_{x\to 0+} \hat{L}^x_t\geq \hat{L}_t\qquad \text{a.s.}
\end{equation}
On the other hand, by the Fubini-Tonelli theorem and the occupation density formula of Proposition \ref{P1}, we see from \eqref{E:Besdens} that for every $x>0$,
$$
\E(\hat L^x_t) = \frac{\alpha  2^{\alpha}}{\Gamma(1-\alpha)}\int_0^t\left( \frac{s}{1-2p}\right)^{\alpha-1}
\exp\left (-\frac{(1-2p)x^2}{2s}\right ) \dd s.
$$
Letting $x\to 0+$, we obtain
\begin{equation}\label{E:expect}
\lim_{x\to 0+}  \E(\hat{L}^x_t) =  \frac{2^{\alpha } }{\Gamma(1-\alpha)  (1-2p)^{\alpha -1}} t^{\alpha}= \E(\hat{L}_t),
\end{equation}
where the second equality is from Theorem \ref{T2}. 
 The statement derives now from \eqref{E:liminf} and \eqref{E:expect}, by  the classical argument of Scheff\'e. 
\end{proof}
One can use Propositions \ref{P1} and \ref{P2} in combination with the scaling property and get an alternative proof of Theorem \ref{T1}(i).

\bibliography{RBES.bib}

\end{document}